\documentclass[12pt,a4paper]{article}

\usepackage{amsmath}
\usepackage{amsthm}
\usepackage{amsfonts}
\usepackage{amssymb}
\usepackage{mathrsfs}
\usepackage{graphicx}
\usepackage{subfigure}

\parskip \smallskipamount

\newtheorem{theorem}{Theorem}[section]
\newtheorem{lemma}{Lemma}[section]

\theoremstyle{remark}

\def\R{\mathbb{R}}

\def\N{\mathbb{N}}

\def\P{\mathbf{P}}
\def\E{\mathbf{E}}

\def\LL{\mathcal{L}}
\def\SS{\mathcal{S}}

\def\HH{\mathcal{H}}

\def\(({\left(}
\def\)){\right)}
\def\as{\textrm{ as }}
\def\forall{\textrm{ for all }}

\renewcommand{\phi}{\varphi}
\renewcommand{\epsilon}{\varepsilon}

\newcommand{\beqr}{\begin{eqnarray*}}
\newcommand{\eeqr}{\end{eqnarray*}}
\newcommand{\bal}{\begin{align*}}
\newcommand{\enal}{\end{align*}}
\renewcommand{\limsup}{\varlimsup}

\newcommand{\oo}[1]{\overline{#1}}
\newcommand{\defn}{:=}

\title{On Upper Bounds for the Tail Distribution of Geometric Sums of Subexponential Random Variables}
\author{Andrew Richards}
\date{}

\begin{document}

\maketitle



\begin{center}
\textit{Department of Actuarial Mathematics and Statistics and the Maxwell Institute for Mathematical Sciences,Heriot-Watt University,Edinburgh}\\
E-mail awr2@hw.ac.uk
\end{center}

\begin{quotation}\small

The approach used by Kalashnikov and Tsitsiashvili for constructing upper bounds for the tail distribution of a geometric sum with subexponential summands is reconsidered.  By expressing the problem in a more probabilistic light, several improvements and one correction are made, which enables the constructed bound to be significantly tighter.  Several examples are given, showing how to implement the theoretical result.

\end{quotation}


\section{Introduction}
This paper was motivated by a series of papers by Kalashnikov and Tsitsiashvili \cite{KT1, KT2}.  In these papers they proposed a method for constructing bounds for the relative accuracy of the asymptotic approximation to the tail distribution of a geometric sum.  In trying to understand their method we found an error in their derivation which led to an incorrect formula (Lemma 1 in \cite{KT1} and Lemma 6 in \cite{KT2}).  However, their general methodology was sound.  We used a more probabilistic approach to deriving equivalent corrected results.  This has simplified the picture, and the more natural probabilistic interpretation has enabled us to achieve much tighter bounds over a greater domain of application.

Let $X_1, X_2,\ldots$ be a sequence of non-negative independent and identically distributed random variables with unbounded support on the positive half-line.  Let $\nu$ be an independent counting random variable with geometric distribution, \[\P(\nu=k)=p(1-p)^{k-1}, \ k \geq 1,\  0<p<1.\]
Then $S_{\nu}=\sum_{i=1}^{\nu}X_i$ is a \emph{Geometric Sum}.  We define $S_0=0.$

Many situations can be modelled as a geometric sum.  Applications include risk processes, ruin theory, queueing theory and reliability models.  The following list of references contains useful surveys: Abate, Choudhury and Whitt \cite{Ab1}, Asmussen \cite{As1}, Cocozza-Thivent and Kalashnikov \cite{CK1}, Feller \cite{Fel1}, Gnedenko and Korolev \cite{Gn1}, Gnedenko and Ushakov \cite{Gn2}, Grandell \cite{Gr1}, Kalashnikov \cite{Kal1, Kal2}.  

In queueing theory, there are many popular models where geometric sums arise naturally.  The first classical example is for a GI/GI/1 queue with stationary waiting time.  It is well known (see \cite{Fel1}) that the waiting time coincides in distribution with the supremum of a random walk with negative drift, $M=\sup_n \sum_{i=1}^n \xi_i$, where the $\xi_i$ are functions of the interarrival times and service times, and have common distribution function $G(x)$.  The supremum of a random walk can then be modelled as a geometric sum $\P(M>x) = (1-p)\P(S_{\nu}>x)$, where $S_{\nu}=\sum_{i=1}^{\nu}\psi_i$ and $p =\P(M=0)$.  The random variables $\psi_i$ appearing in the geometric sum have the distribution of the first overshoot of the random walk over level zero, conditional on such an overshoot occurring.  The distribution of the $\psi_i$ does not coincide with $G$.  However we can obtain bounds on the distribution of the $\psi_i$, (see, for example, Borovkov \cite {Bor1}).  In \cite{Bor1}, Chapter 4, Section 22, Theorem 10 we find
\begin{equation}\label{qu}
c_1\oo{G^I}(x) \leq \P(\psi_i >x) \leq c_1\oo{G^I}(x) + c_2\oo{G}(x),
\end{equation}
where $\oo{G(x)} = 1 - G(x)$ and $\oo{G^I}(x) = \frac{1}{\E(\xi_1)} \int_x^{\infty} \oo{G}(y)dy$, and $c_1$ and $c_2$ may be estimated under further assumptions.  We shall give an example of this in section 4. 

The tail distribution of the geometric sum is the object of interest.  In the light-tailed case the situation is straightforward and an upper bound for the tail distribution is well known, the Cram\'er Lundberg upper bound.  However, when the summands have a subexponential distribution it is much more difficult to find an appropriate upper bound for the distribution.   Analytical techniques, simulation and recursive techniques have been used to estimate the tail.  The tail is usually approximated by its asymptotic form (see, for example, Embrechts and Veraverbeke \cite{Em1}), but this approximation can be very poor, and therefore it is of use to have an upper bound for the relative accuracy of the asymptotic approximation.  A complete list of references for works on the bounds of the tail of the geometric sum distribution would be too long, but among important works are: Asmussen, Binswanger and Hojgaard \cite{As2}, Dickson and Waters \cite{DW1}, Kalashnikov \cite{Kal1,Kal2}, Lin \cite{Lin1}, Willmot \cite{Wil1}, Willmot and Lin \cite{WL1}.

The rest of this paper is structured as follows.  In section 2 we lay the ground work and derive an auxiliary result concerning subexponential distributions.  In section 3 we obtain the required upper bound for the relative accuracy of the asymptotic expression.  In section 4 we apply the result to a number of subexponential distributions.

\section{Some Preliminary Results}

We consider distributions $F(x)$ with unbounded support on the positive half-line.  We will write the tail distribution as $\oo{F}(x)=1-F(x)$, and $F^{*n}(x)$ as the $n$-fold convolution.

We recall the following standard definitions and classes of distributions.  For any functions $a(x) \geq 0, b(x) \geq 0$,
\begin{enumerate}
\item[i)] $a(x)=o(b(x))$ means $\lim_{x \to \infty} \frac{a(x)}{b(x)} = 0$;\\
\item[ii)] $a(x)=O(b(x))$ means $0 < \limsup_{x \to \infty} \frac{a(x)}{b(x)} < \infty$;\\
\item[iii)] $a(x) \sim b(x)$ means $\lim_{x \to \infty} \frac{a(x)}{b(x)} = 1$.
\end{enumerate}
The class $\HH$ of heavy-tailed distributions consists of those distributions $F$ for which $\int_0^{\infty} e^{tx}F(dx) = \infty$ for all $t>0$.
A subclass of the heavy-tailed distributions is the class $\LL$ of long-tailed distributions, defined by the requirement that
\[\lim_{x \to \infty} \frac{\oo{F}(x-1)}{\oo{F}(x)}=1.\]
It is clear that, for any $a \in \R$, $\oo{F}(x-a) \sim \oo{F}(x)$.
This in turn implies the existence of an increasing concave function $h(x)<x/2$ such that $h(x) \to \infty$ as $x \to \infty$, see Foss and Zachary \cite{FZ}, and
\begin{equation}\label{h}
\oo{F}(x-h(x)) \sim \oo{F}(x).
\end{equation}
The class $\SS$ of subexponential distributions was introduced by Chistyakov in 1964 \cite{Chist}.  A distribution $\oo{F}$ on the positive half-line belongs to $\SS$ if and only if
\[\oo{F^{*2}}(x) \sim 2 \oo{F}(x).\]
The class $\SS$ is a subclass of $\LL$ (see Kl\"{u}ppelberg \cite{Klup}).
By induction we also have $\oo{F^{*n}}(x) \sim n \oo{F}(x)$ for any $n \in \N$.
If $F \in \SS$, then for any function $g(x) \to \infty$ as $x \to \infty$, $g(x) \leq x$,
\begin{equation}\label{intse}
\int_{g(x)}^{x} \frac{\oo{F}(x-y)}{\oo{F(x)}}F(dy) \to 0 \as x \to \infty.
\end{equation}
This result, to be found in \cite{FKZ}, follows from the observation that \eqref{intse} is valid for any $g(x) = h(x)$, where $h(x)$ satisfies \eqref{h}, and if \eqref{intse} is valid for some $g(x)$, then it is valid for any $g_1(x)$ such that $g(x) \leq g_1(x) \leq x$.
For convenience we define the following notation.
\begin{align}
K_F(x,r) \equiv K(x,r) & \defn \frac{\oo{F}(x-r)}{\oo{F}(x)}-1,\\
J_F(x,r) \equiv J(x,r) & \defn \int_r^{x-r}\frac{\oo{F}(x-y)}{\oo{F}(x)}F(dy).
\end{align}
Both these expressions have natural interpretations for, respectively, long-tailed and subexponential distributions.

We note that, for $F \in \SS$,
\begin{enumerate}
\item[i) ] $J(x,r) \to \oo{F}(r)$ as $  x \to \infty$ for fixed $r$;
\item[ii) ] $J(x,r) $ is monotonically decreasing in $r$ for fixed $x$;
\item[iii) ] $K(x,r) \to 0 \as x \to \infty$ for fixed $r$;
\item[iv) ] $K(x,r)$ is monotonically increasing in $r$ for fixed $x$.
\end{enumerate}

\begin{lemma}\label{Lem1}
Let $F \in \LL$. Then $F \in \SS$ if and only if for any $c>1$ there exists a function $h(x) \equiv h(x,c)$ and a constant $x_1 \equiv x_1(c) >0$ such that, for $x \geq x_1$,
\begin{enumerate}
\item[i) ] $h(x)$ is an increasing concave function ;
\item[ii) ] $h(x) \to \infty$ as $x \to \infty$;
\item[iii) ] there exists $x_0 \geq 0$ such that $h(x) \leq x/2$ for $x>x_0$;
\item [iii) ] $J(x,h(x)) \leq c \oo{F}(h(x))$;
\item[iv) ] $K(x,h(x)) \leq \oo{F}(h(x))$.
\end{enumerate}
\end{lemma}

\begin{proof}
Let $F \in \SS$ and $c>1$.  Then there exists $k \in \N$ such that
\[k=\min\((j \in \N : \inf_{x \geq j} \frac{\oo{F}(x+1)}{\oo{F}(x)} > \frac{1}{c}\)).\]
Let $x_{-1}=x_0=0$.

For fixed $t \geq k$, $J(x,t) \to \oo{F}(t) < c\oo{F}(t+1)$.  So, for $r \in \N$, we can inductively define an unbounded, increasing sequence:
\[x_r=\inf(x: x-x_{r-1} \geq x_{r-1}-x_{r-2}\geq 3, \sup_{y \geq x} J(y, k+r-1) \leq c\oo{F}(k+r),
	\sup_{y \geq x} K(y, k+r) \leq \oo{F}(k+r)).\]
Now define
\[h(x) = \left\{ \begin{array}{cc} kx/x_1 & x<x_1,\\ k+ r - 1 + \frac{x-x_r}{x_{r+1}-x_r} &  x \in [x_r, x_{r+1}), r \in \N.
\end{array} \right. \]
For $x \in [x_r, x_{r+1}), r \in \N,$ we have
\[k+r-1 \leq h(x) < k+r.\]
So,
\[J(x,h(x)) \leq J(x,k+r-1) \leq c\oo{F}(k+r) \leq c\oo{F}(h(x))\]
and
\[K(x,h(x)) \leq K(x,k+r) \leq \oo{F}(k+r) \leq \oo{F}(h(x)).\]
By construction increasing, unbounded above and concave, and, for $x>x_1+6K$, $h(x)\leq x/2$.

On the other hand, if $F$ is long-tailed, and such a function $h(x)$ exists, then by Proposition 2 of \cite{AFK} $F$ is subexponential.
\end{proof}

\section{Upper Bound for the Relative Accuracy}

The main goal of this section is to obtain Theorem \ref{Th1}.  Our intention was to make the statement and proof of this theorem look very similar to Theorem 3 in \cite{KT2}.  However, the quantities defined in the statement of our Theorem \ref{Th1} are not identical to the corresponding quantities in Kalashnikov and Tsistsiashvili's paper.

The structure of this section is as follows.  First we discuss the accuracy of the asymptotic approximation \eqref{1}.  We construct an upper bound for $\P(S_{\nu}>x|\nu>1)$.  We introduce some notation for two common expressions associated with subexponential distributions, and use these to construct an upper bound for the relative accuracy $\Delta(x)$.  We then introduce our test function $g(x)$, whose existence is guaranteed by Lemma \ref{Lem1} and re-express the upper bound for $\Delta(x)$ in terms of $C(x)$, the ratio between $\Delta(x)$ and the test function $g(x)$.  Finally we prove Theorem \ref{Th1}.

\subsection{Accuracy of the Asymptotic Approximation}

It is well known (see, for example, \cite{Em1}) that asymptotically
\begin{equation}\label{1}
\P(S_{\nu}>x) \sim \E\nu \oo{F}(x).
\end{equation}

In general the accuracy of the asymptotic relation \eqref{1} is poor for small values of $x$.  For instance, if we consider a Pareto(5,1) distribution, and take $p=0.2$, so
\[\oo{F}(x) = \left\{
\begin{array}{cc} 1 &  \textrm{for  } x \leq 1,\\ x^{-5} & \textrm{for  } x > 1 \end{array} \right.\]
then, estimating the tail of the geometric sum using naive simulation,
\begin{align*}
\P(S_{\nu}>30) & = 0.00547 \\
\E(\nu)\P(X>30)& = 2.06 \times 10^{-7}.
\end{align*}
This gives us a relative error of more than $26000$ in using the asymptotic expression.  Hence the need to establish an upper bound for the relative accuracy of the asymptotic relation.

Let
\[\Delta_F(x) \equiv \Delta(x) = \frac {\P(S_{\nu} >x) - \E\nu \oo{F}(x)}{\E\nu \oo{F}(x)}.\]
and $\Delta[a,b] = \sup_{a \leq x \leq b} \Delta(x).$

\subsection{Upper bound for $\P(S_{\nu}>x|\nu>1)$}

We are now ready to construct an upper bound for $\P(S_{\nu}>x|\nu>1)$.

Choose $h(x)$, such that $\oo{F}(x-h(x)) \sim \oo{F}(x)$.
We consider the total probability formula:
\begin{equation}\label{2}
\P(S_{\nu} >x)=\P(X_1>x)\P(\nu=1)+\P(S_{\nu}>x|\nu>1)\P(\nu>1).
\end{equation}

Conditional on $\nu > 1$, $S_{\nu} - X_1$ has the same distribution as $S_{\nu}$.  Let $\widetilde{S}_{\nu}$ be a random variable, independent of $\{X_i\}$ with the same distribution as $S_{\nu}$. 
Denote the event $\{ S_{\nu} >x,\nu >1 \} \defn B(x)$.
We can then partition the event $\{S_{\nu}>x,\nu>1\}$ as follows:
\begin{align*}
B(x) = \{X_1 \leq h(x), B(x)\} &\cup \{h(x) < X_1 \leq x-h(x), B(x)\} \\
&\cup \{X_1>x-h(x), B(x)\}.
\end{align*}
We then have
\begin{align*}
	\P(S_{\nu} >x|\nu >1) \leq & \P(X_1 \leq h(x),\widetilde{S}_{\nu}> x-h(x))\\
	& + \P(h(x) < X_1 \leq x-h(x), \widetilde{S}_{\nu} > x - X_1)  + \P(X_1 > x-h(x)). 
\end{align*}
We know that
\[\P(S_{\nu} > x) = \frac{1}{p}(1+ \Delta(x))\oo{F}(x).\]
Hence,
\begin{align*}
	\P(S_{\nu} >x|\nu >1)  \leq & \P(X_1 \leq h(x))\P(\widetilde{S}_{\widetilde{\nu}} > x-h(x))\\
	 & + \int_{h(x)}^{x-h(x)}\P(\widetilde{S}_{\widetilde{\nu}}>x-y)\P(X_1 \in dy) + \P(X_1 > x-h(x))\\
	\leq & \frac{1}{p}(1+ \Delta(x-h(x)))\oo{F}(x-h(x))F(h(x))\\
	& + \frac{1}{p}(1+ \Delta[h(x),x-h(x)]) \int_{h(x)}^{x-h(x)} \oo{F}(x-y)F(dy) + \oo{F}(x-h(x)).
	\end{align*}

\subsection{An Upper Bound for $\Delta(x)$}

To simplify our expressions we recall the quantities defined in (3) and (4), and denote $J \equiv J(x,h(x))$ and $K \equiv K(x,h(x))$.
We note that, since $F \in \SS$, we may choose a concave function $h(x) \to \infty$ such that both $J$ and $K$ converge to $0$ as $x \to \infty$.

Using such a function $H$ we can now construct an upper bound for the relative accuracy $\Delta(x)$.

	We have
	\begin{align*}
		\frac{\P(S_{\nu}>x)}{\oo{F}(x)} \leq & \frac{1-p}{p}\(( (1+\Delta(x-h(x)))(K+1)F(h(x)) + (1+\Delta[h(x),x-h(x)])J \))\\
		& +(1-p)(K+1)  + p,
	\end{align*}

and so,
\begin{align}
	\Delta(x) \leq & (1-p)\((\Delta(x-h(x))(K+1) F(h(x)) + \Delta[h(x),x-h(x)]J \))\nonumber\\
	& + (1-p)J + (1-p^2)K - (1-p)(K+1)\oo{F}(h(x)).
	\end{align}

Our upper bound for $\Delta(x)$ should be a monotonic function $\Delta_u(x)$, decreasing to 0, such that, for some constant $b$,
\[\Delta(x) \leq \Delta_u(x) \quad \forall x \geq b.\]

\subsection{The Test Function}

As in \cite{KT2} we now introduce our test function $g(x)$. 

We require $g(x)$, which will depend on our choice of $h(x)$, to have the following properties as $x \to \infty$:
\begin{align}
g(x) & \to 0 \quad \textrm{monotonically};\\
\max(J(x,h(x)),K(x,h(x))) & = O(g(x));\\
g(x-h(x)) & \sim g(x);
\end{align}
The existence of such a function is again guaranteed by Lemma \ref{Lem1} since we may choose $g(x) \defn \oo{F}(h(x))$, with $h(x)$ concave as in Lemma \ref{Lem1}.  The relation $h(x-h(x)) \geq h(x) - h(h(x))$, for large enough $x$, then implies (10).
Define
\[
C(x) = \max\((0,\frac{\Delta(x)}{g(x)}\)) \quad \mbox{and} \quad C[a,b] =\max_{a \leq x \leq b}C(x).\]

Now, if we set
\begin{align*}
f_1(x) & = \frac{(1-p)g(x-h(x))(K+1)F(h(x))}{g(x)},\\
f_2(x) & = \frac{(1-p)g(h(x))J}{g(x)},\\
f_3(x) & = \frac{(1-p)J + (1-p^2)K - (1-p)(K+1)\oo{F}(h(x))}{g(x)}
\end{align*}
we can rewrite (7) as
\begin{equation}\label{4}
C(x) \leq (f_1(x)+f_2(x))C[h(x),x] + f_3(x).
\end{equation}
Also, $f_1(x) \to 1-p$, $f_2(x) \to 0$, and $f_3(x)$ is bounded from above.  Let $\delta \equiv \delta(x)=\sup_{y \geq x}((f_1(y)+f_2(y))$.
Hence we can find $b>0$ such that 
\begin{equation}\label{5}
(f_1(x)+f_2(x)) \leq \delta(b)  < 1 \textrm{ and } f_3(x) \leq \phi(b) \equiv \phi \forall x \geq b.
\end{equation}

\subsection{The Key Result}

Following the arguments in \cite{KT2} we can now prove:
\begin{theorem}\label{Th1}
Let $F \in \SS$, $h(x)$ satisfy condition \eqref{h}, and $g(x)$ satisfy conditions (8), (9) and (10).  Then, there exists $b>0$, $0<\delta<1$ and $\phi > 0$, such that for all $x \geq b$,
\[\Delta(x) \leq \Delta_u(x) \defn Cg(x)\]
where
\begin{equation}\label{6}
C = \max\left( \frac{\phi}{1-\delta},\phi+\delta C[h(b),b]\right).
\end{equation}
\end{theorem}
\begin{proof}
For $x \geq b$ it is clear that $\Delta(x) \leq C(x)g(x) \leq C[b,x]g(x).$\\
From \eqref{4} we know that $C(x) \leq (f_1(x)+f_2(x))C[h(x),x] + f_3(x).$  Hence, with $b, \delta$ and $\phi$ as in \eqref{5},
\[C[b,x] \leq \delta C[h(b),x] + \phi.\]
Now, $C[h(b),x] = \max(C[h(b),b],C[b,x])$.  

If $C[h(b),x] = C[h(b),b]$ then
\[C[b,x] \leq \delta C[h(b),b] +\phi.\]
If $C[h(b),x] = C[b,x]$ then
\[C[b,x] \leq \frac{\phi}{1-\delta}.\]
This completes the proof.
\end{proof}

This result can be applied to the queueing problem described in the introduction by using the estimate \eqref{qu}. We have upper and lower bounds for the tail distribution and  we call these $\oo{F^+}$ and $\oo{F^-}$.  We estimate the relative error, $\Delta_{F^+}$, in using the asymptotic approximation $\P(S_{\nu}>x) \approx \E(\nu)\oo{F^+}(x)$.  All the results for Theorem 1 follow through using $F^+$ in place of $F$.  We need to evaluate $C(x)$ for $h(b) \leq x \leq b$, and we propose to do this using naive simulation.  We shall bound $C(x)$ above by considering a geometric sum random variable $S_{\nu}^+$ whose increments have distribution $F^+$, so that
\[C(x) \leq C^+(x) \defn \((p\frac{\P(S_{\nu}^+ > x)}{\oo{F^+}(x)} -1\))/g(x).\]

\section{Applying the Result}

Some numerical estimation of the tail distribution of the geometric sum must be done in order to evaluate $C[h(b),b]$.  The greater the value of $b$, the more accurate the upper bound becomes.  However, this comes at the greater computational cost of numerically evaluating the tail of the distribution.  A compromise has to be struck between the tightness of the upper bound and the resources one is willing to invest in evaluating the tail.

A critical part of the procedure is the choice of the test function $g(x)$, which itself depends on the choice of $h(x)$.  In \cite {KT2} the function $g(x)$ was chosen as a function in closed form over the whole range of values of its argument.  We observe, however, that in evaluating $C[h(b),b]$, we know the (numerically) exact value of $\Delta(x)$ in the range $h(b) \leq x \leq b$, and we wish to use this information in our choice of $g(x)$ when appropriate.  Our strategy in applying the result is therefore as follows.  
\begin{enumerate}
\item[1) ] Decide what resources are available for estimating the tail distribution to a suitable degree of accuracy.  Given the available resources, define $B$ as the maximum value for which we estimate $\Delta(B)$ and numerically evaluate $C[h(B),B]$.
\item[2) ] Determine the class of functions that will do for $h(x)$.
\item[3) ] Estimate $J(x,h(x))$ and $K(x,h(x))$ in the range $h(B) \leq x \leq B$.
\item[4) ] Choose monotonically decreasing $g(x)$, which will depend on our particular choice for $h(x)$, and which may incorporate our numerical knowledge of $C(x)$, such that 
$g(x) = O( \max(J(x,h(x)),K(x,h(x))))$.
\item[5) ] If $\sup_{x \geq B}\delta(x) <1$, we take $b=B$, and find the corresponding value of $\phi$.  If $\sup_{x \geq B}\delta(x) \geq 1$ either the procedure has failed, or we must be prepared to use a larger value of $B$.
\item[6) ] Calculate $C$.
\end{enumerate}
Some comments on these steps will be useful.  

In step 2 when choosing $h(x)$, there is a tension involved between the relative rates of decay of $J(x,h(x))$ and $K(x,h(x))$.  The larger $h(x)$ is the smaller $J(x,h(x))$ becomes, but the larger $K(x,h(x))$ becomes, and vice versa.  We can change the rates of decay of $J(x,h(x))$ and $K(x,h(x))$ by scaling $h(x)$ by some numerical factor without affecting the asymptotic decay rate of $g(x)$.  

In Step 4 we will generally want to choose $g(x)$ in order to cause the upper bound for the relative accuracy to decay to zero as fast as possible.  We may also want to incorporate the information we have already calculated for $C[h(B),B]$.   The fastest asymptotic decay rate for $g(x)$ is obtained by an optimal choice of $h(x)$.  However, we can use the information we have gathered in calculating the numerically exact value of $\Delta(x)$ in the range $[h(b),b]$ by constructing a monotonically decreasing version of it, $\Delta_m(x) \defn \sup_{x \leq y \leq b^*} \Delta(y)$,  for $x<b^*$, where $b^*<b$ is chosen to minimize the value of $C$.  Thus, once we know the optimal asymptotic function $g(x)$, we instead use $g_1(x)$:
\[ g_1(x)= \left\{
	\begin{array}{cc}
			\Delta_m(x) & \textrm{for  } x<b^*,\\
			Kg(x) & \textrm{for  }  x \geq b^*,
	\end{array} 
\right.
\]
for some constant $K$ chosen to make $g_1(x)$ continuous at $b^*$.

Hence we have a two parameters that we can adjust, the scale factor for $h(x)$, which alters the balance between the decay of $J(x,h(x))$ and $K(x,h(x))$; and the value of $b^*$, which allow us, given our chosen value of $B$, to minimize $C$, and hence tighten the upper bound.

We will now show how to apply our result to Pareto and Weibull distributions with various parameters.
The values of $J(x,h(x)$ were estimated using the integrate function in R, and the naive simulation to estimate $C(x)$ was also performed in R using samples of size $5 \times 10^6$.  The bounds we calculate (labelled KT bounds in the graphical displays to reflect the original source of this method in the work of Kalashnikov and Tsitsiashvili) are compared to values of the relative error that were calculated using a discretized Panjer algorithm with bandwidth of 0.005.

\subsection{Pareto Distribution}
We will consider Pareto Distributions of the following form:
\[ F(x)= \left\{
	\begin{array}{cc}
			0 &  \textrm{for  } x<1\\
			1- x^{-\alpha} &  \textrm{for  } x \geq 1
	\end{array} 
\right.
\]
where $\alpha > 1$.  We now follow the steps above.  
The choice of $h(x)$ is determined by the requirement that $K(x,h(x)) \to 0$. This occurs if and only if $h(x) = o(x)$.  We then have
\begin{align*}
K(x,h(x)) & = \((1-\frac{h(x)}{x}\))^{-\alpha}-1   \leq \frac{\alpha h(x) x^{\alpha}}{(x-h(x))^{\alpha+1}},\\
J(x,h(x)) & \leq 2\int_{h(x)}^{x/2} \frac{\oo{F}(x-y)}{\oo{F}(x)}F(dy) + \frac{\oo{F}^2(x/2)}{\oo{F}(x)}	  \leq 2\((\frac{2}{h(x)}\))^{\alpha} + \((\frac{4}{x}\))^{\alpha}.
\end{align*}
The simplest form of $h(x)$ is $h(x) = x^{\beta}$ for $0<\beta<1$.  We can then choose 
\[g(x) = x^{-\min(\alpha \beta, 1- \beta)},\]
which ensures that $g(x) = O( \max(J(x,h(x)),K(x,h(x)))).$
If we want to make $g(x)$ decay as fast as possible the optimal choice for $h(x)$ will have $\beta = 1/(1+\alpha)$.
\subsubsection{Example 1}
We shall take $\alpha=2.2$ and $p=0.5$.  We use a discretized Panjer recursion to estimate the tail of the distribution, and assume that our resources allow us to estimate this up to $B=100$ using a bandwidth of $0.005$.  If we follow the approach of \cite{KT2} we take $h(x)=x^{1/3.2}$ and 
\[g(x)=x^{-2.2/3.2}=x^{-0.6875}.\]
We find that $\delta(100) = 0.786$, and $C[5,100]=14.4$.  This results in $C=13.2$, giving 
\[\Delta(x) \leq 13.2x^{-0.6875}, \quad x>100,\]
as the (corrected) Kalashnikov-Tsitsiashvili bound.  A small improvement can be made by taking $h(x)=1.70x^{1/3.2}$, with a consequent value of $C=12.9$.  A further improvement can be made, at some cost to computational time , by choosing the test function $g(x)$ to be equal to $\Delta_m(x)$, the monotonically decreasing version of the exact value of $\Delta(x)$, up to some value $b^* \leq B$, and then ensuring continuity at $b^*$.  The best that can be obtained with our value of $B=100$ is to take $h(x)=1.14x^{1/3.2}$, and $b^*=21.3$, so that
\[g(x) = \left\{	\begin{array}{cc}  \Delta_m(x) & \textrm{for  } x\leq 21.3,\\
	8.52 x^{-0.6875} &  \textrm{for  } x>21.3. \end{array} \right.\]
 This results in the upper bound
\[\Delta(x) \leq 8.53 x^{-0.6875}, \quad x>100.\]
  The logarithm (base 10) of this upper bound has been graphed  in Figure \ref{fig:a} (labelled as Upper Bound), along with the logarithm of the numerically exact result obtained from the Panjer recursion (labelled Exact).

\subsubsection{Example 2}
Now consider $\alpha=2.2$ and $p=0.2$.  Once again we take $B=100$.  If we follow the methodology of \cite{KT2} and take $h(x)=x^{1/3.2}$ and $g(x) = x^{-2.2/3.2}$, we find that $\min(n \in \N: \delta(n)<1)=1085$ which not only is greater than our chosen $B$, but impracticably large in any event. No improvement can be made just by scaling $h(x)$.  However, if we take $h(x)=1.054x^{1/3.2}$ and adjust $g(x)$ to coincide with $\Delta_m(x)$ for $x<b^*=27.1$, we arrive at
\[g(x) = \left\{	\begin{array}{cc}  \Delta_m(x) &  \textrm{for  } x\leq 27.1,\\
	126 x^{-0.6875} &  \textrm{for  } x>27.1. \end{array} \right.\]
This gives the upper bound
\[\Delta(x) \leq 179.85 x^{-0.687}, \quad x>100,\]
which is shown in Figure \ref{fig:b}.

\subsubsection{Example 3}
\begin{figure}[tb]
     \begin{center}
     \subfigure[Ex 1. Pareto, $\alpha = -2.2$, 
     $p=0.5$,  \quad $\Delta_u(x)=8.53 x^{-0.6875}$]{
          \label{fig:a}
                 \includegraphics[width=6cm,height=6cm]{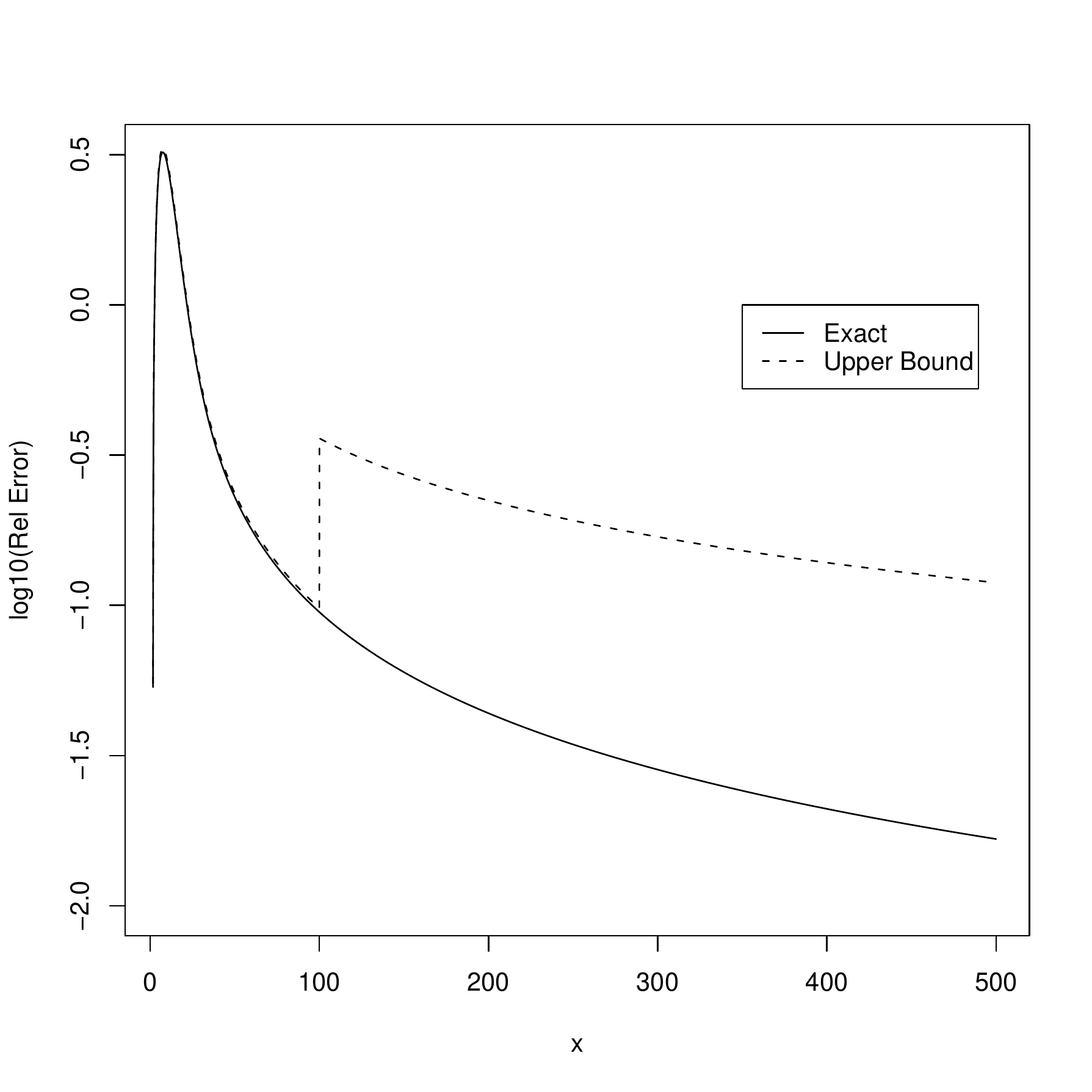}}
     \hspace{0.5cm}
     \subfigure[Ex 2. Pareto, $\alpha = -2.2$, 
     $p=0.2$, \quad $\Delta_u(x)=179.85 x^{-0.6875}$]{
          \label{fig:b}
                \includegraphics[width=6cm,height=6cm]{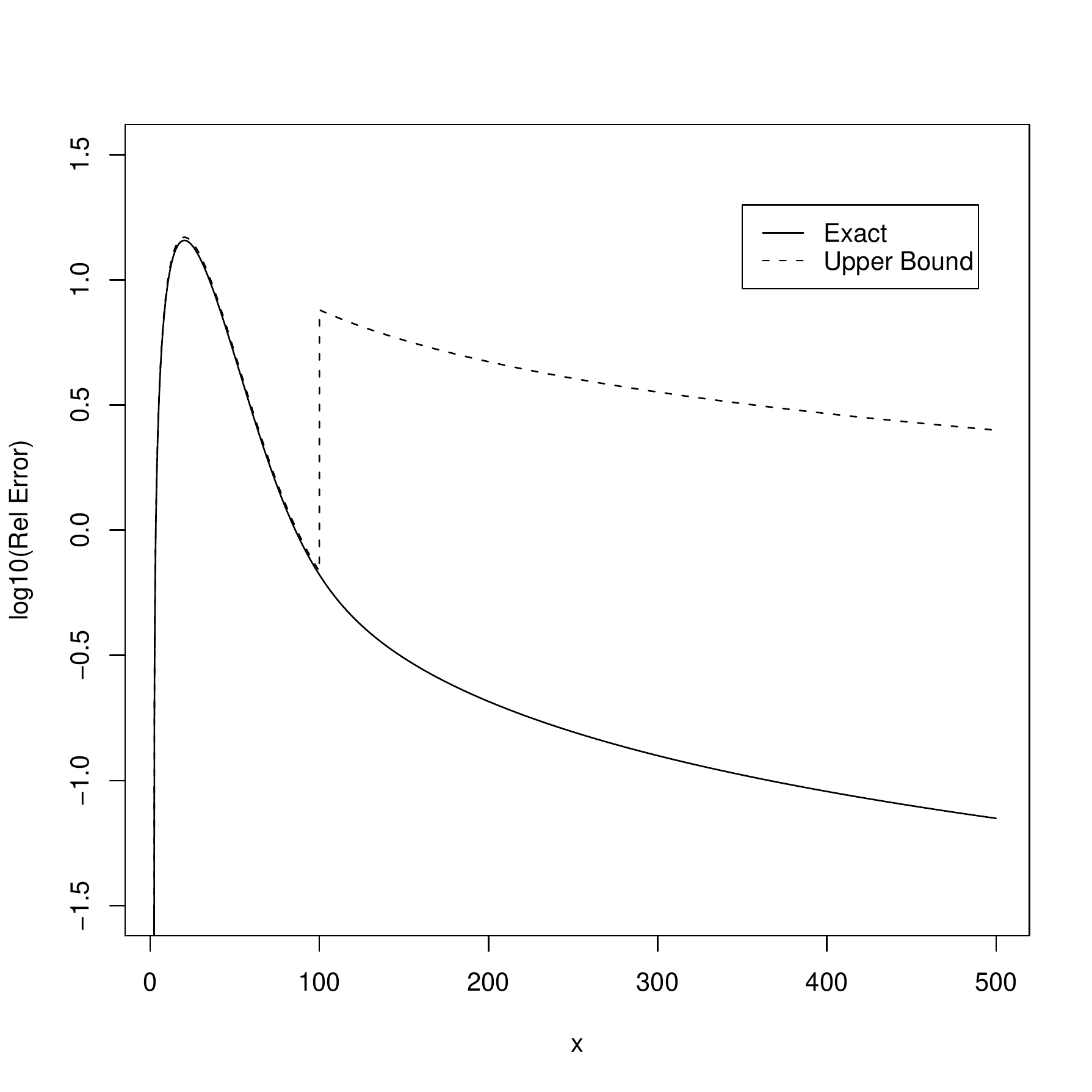}}\\
     \vspace{0.5cm}
     \subfigure[Ex 3. Pareto, $\alpha = -5$,  $p=0.5$,\quad $\Delta_u(x)=93.7x^{-5/6}$]{
           \label{fig:c}
                \includegraphics[width=6cm,height=6cm]{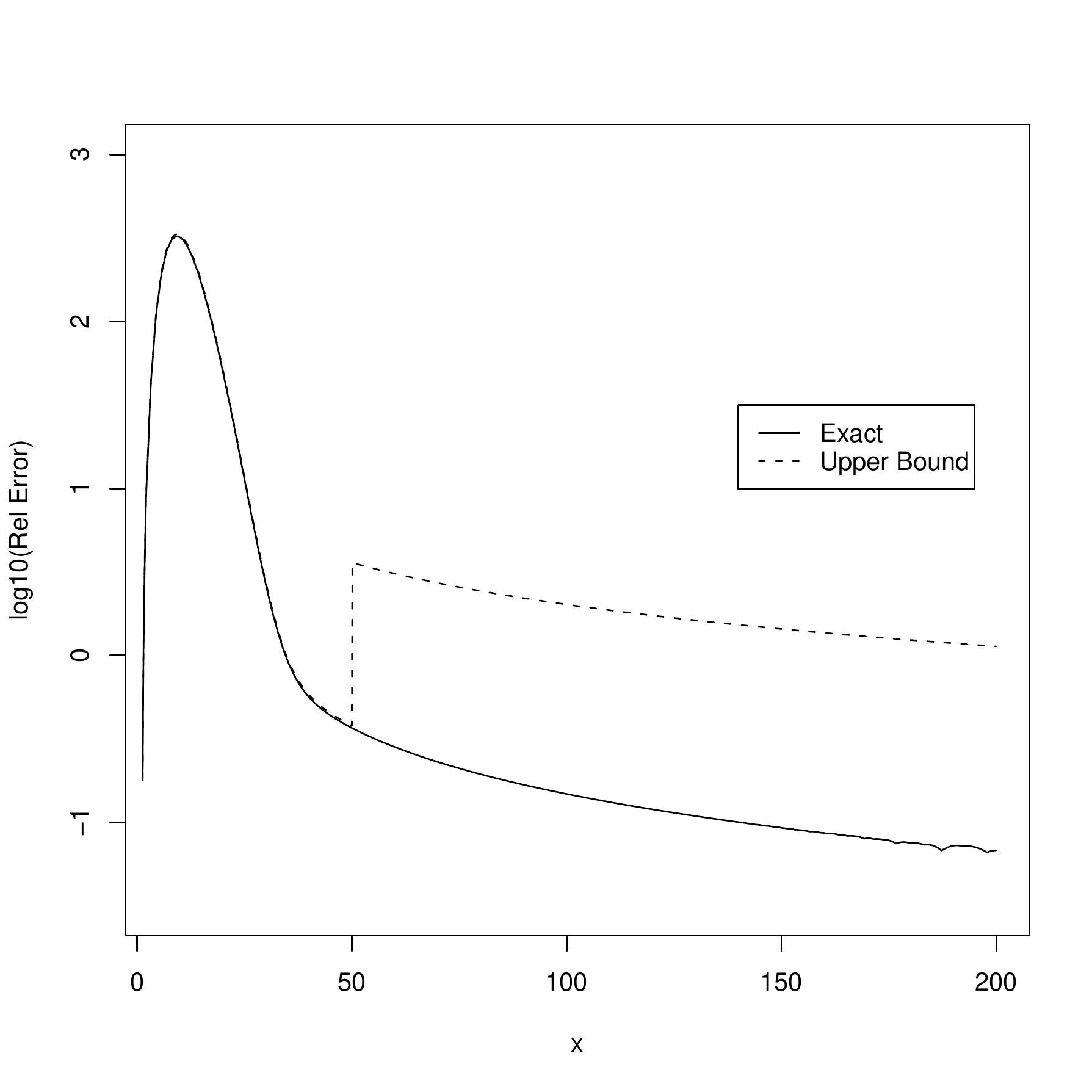}}
                \hspace{0.5cm}
     \subfigure[Ex 5. Weibull, $\beta =0.5$,  $p=0.5$,\quad $\Delta_u(x)= 2.952 K(x,h(x))$]{
           \label{fig:d}
                \includegraphics[width=6cm,height=6cm]{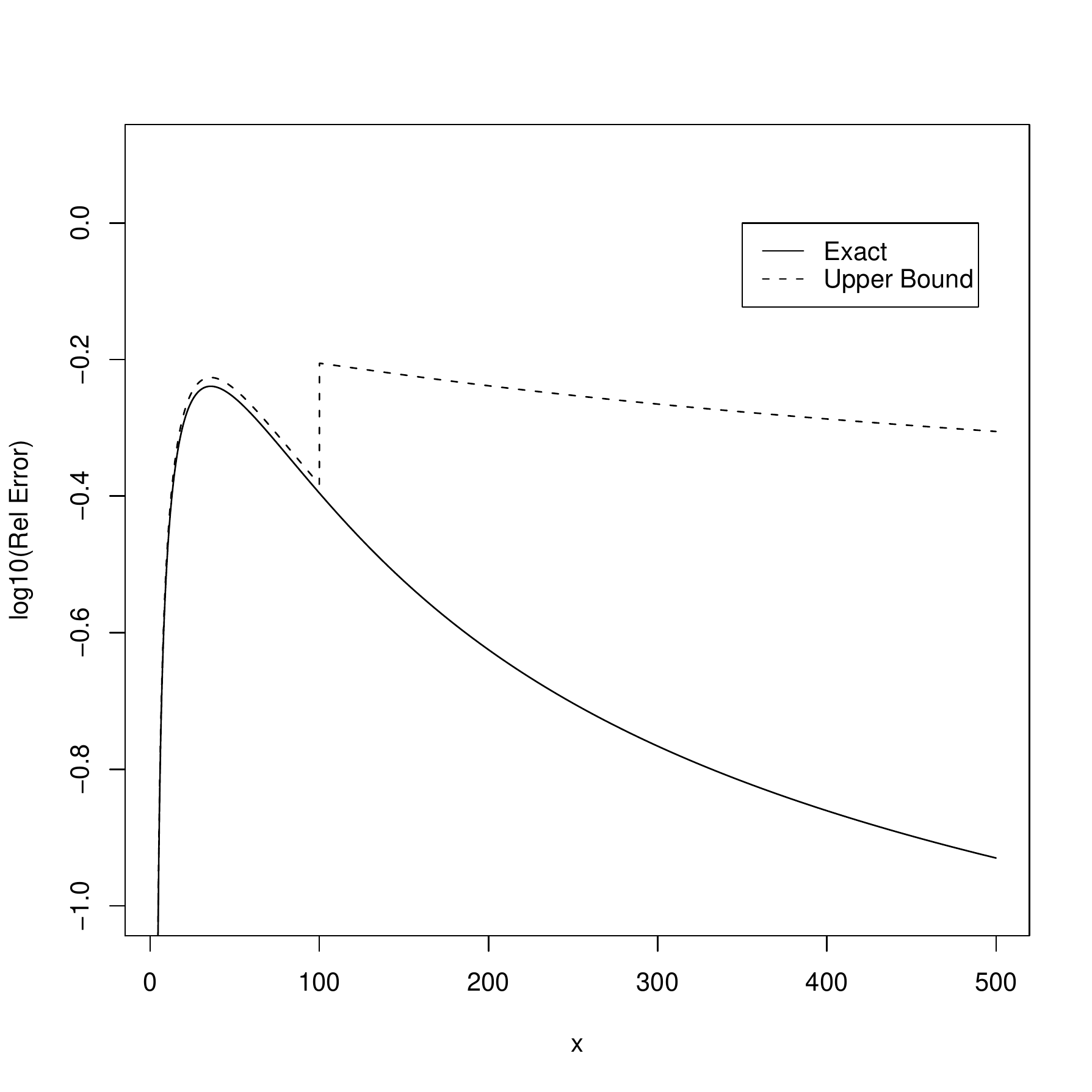}}
     \caption{Examples 1,2,3,5. Plot of $\log_{10}(\Delta(x))$ and $\log_{10}(\Delta_u(x))$.}
     \label{fig:1}
     \end{center}
\end{figure}

When the Pareto distribution in question is lighter tailed the asymptotic approximation becomes drastically less good for moderate values of the argument, and it requires more resources to compute numerically the tail distribution of the geometric sum for higher values of the argument.  As our example we take  $\alpha=5,p=0.5$ and shall again perform this numerical exercise using a discretized Panjer algorithm, but with $B=50$ and a bandwidth of $0.002$.  Proceeding as in \cite{KT2}, we take $h(x)=x^{1/6}$ and $g(x)=x^{-5/6}$.  We find that $\delta(50)=0.996$, and $C=2215$.  This large value of $C$ reflects the fact that the largest errors in the asymptotic expression for the tail of the sum occur very early on (as can be seen from \ref{fig:c}).  By adjusting $h(x)$ to $h(x)=1.46x^{1/6}$ we can reduce $C$ to $C=1662$.  However, in order to avoid the initially large values of $\Delta(x)$ we take $b^*=27.1$, $h(x)=1.94x^{1/6}$ and
\[g(x) = \left\{	\begin{array}{cc}  \Delta_m(x) &  \textrm{for  } x\leq 27.1,\\
	93.93x^{-5/6} &  \textrm{for  } x>27.1, \end{array} \right.\]
which yields a dramatic improvement, giving
\[\Delta(x) \leq 93.7x^{-5/6}, \quad x>50.\]
This bound is shown in Figure \ref{fig:c}.

\subsubsection{Example 4}
For this example we return to the example of the classical GI/GI/1 queue as described in the introduction.  Suppose that we have obtained an upper bound for the distribution of $\psi_i$  as in \eqref{qu} with distribution $F^+$, so that $\P(\psi_i>x) \leq \oo{F^+}(x) = \frac{1}{3}x^{-2} + \frac{2}{3}x^{-3} \textrm{  for  } x>1$, and that we have a value of $p=0.5$.  We estimate an upper bound for the tail of the distribution of the geometric sum using a naive simulation of a geometric sum with increments having distribution $\oo{F^+}$ and decide we have enough resources to $B=80$.  We choose $h(x) = x^{1/3}$ and $g(x)=x^{-2/3}$.  We find that $C=13$.  Hence
\[\Delta(x) \leq 13x^{-2/3}, \quad x >80.\]

\subsection{Weibull Distribution}
We will consider Weibull Distributions of the following form:
\[ F(x)= \left\{
	\begin{array}{cc}
			0 &  \textrm{for  } x\leq 0,\\
			1 - e^{-x^{\beta}} &  \textrm{for  } x>0
	\end{array} 
\right.
\]
where $0<\beta <1$.
Straightforward calculations show that $h(x)$ must be chosen so that $h(x)=o(x^{1-\beta})$.  Further calculations show that, if we take $(\log(x))^{1/\beta}=O(h(x))$, then $J(x,h(x)) = O(K(x,h(x))$, and hence we may take $g(x)=K(x,h(x))$.   

\subsubsection{Example 5}

For this example we take $\beta=0.5, p=0.5.$  The numerical calculations were done using a Panjer recursion with bandwidth $0.002$, and $B$ was taken to be $B=100$.  The optimal function for $h(x)$ is $h(x)=(\log(x))^2$.  Applying the methodology in \cite{KT2}, we find that $\delta(x) >1$ for $x<1660$, and hence no results can be obtained (for $B<1660$).  However, by taking $h(x)=0.179 (\log(x))^2$, we obtain 
\[\Delta(x) \leq 2.952 K(x,h(x)), \quad x\geq 50,\]
where $K(x,h(x))=\exp(\sqrt{x}-\sqrt{x-0.179 (\log(x))^2})-1$.  This Weibull bound is shown in Figure \ref{fig:d}.  Because of the lack of a very sharp peak in $\Delta(x)$, no further improvement can be obtained by incorporating the numerical values of $\Delta(x)$ into the function $g(x)$.

\vspace{1cm}

\bfseries\large{Acknowledgements}

\mdseries\normalsize
I am grateful to Serguei Foss for bringing the original papers by Kalashnikov and Tsitsiashvili to my attention, for the many helpful discussions he had with me, and for his much appreciated advice, and to Onno Boxma for his constructive criticism. Also, I would like to thank the anonymous referee for the careful reading of the drafts of this paper and the many helpful comments.

\end{document}